\nonstopmode
\documentclass[a4paper, 10pt]{amsart}
\usepackage{latexsym}
\usepackage{fancyhdr}
\usepackage{amsmath, amssymb}
\usepackage[ansinew]{inputenc}
\usepackage[all]{xy}
\usepackage{verbatim}
\usepackage{psfrag}
\usepackage{epsfig}
\swapnumbers
\newenvironment{proclaim}[1]{\par\medskip\noindent{\bf #1.}\it}{\par\smallskip}
\newenvironment{demo}[1]{\par\smallskip\noindent{\bf #1.}}{\par\smallskip}

\newcommand{\nmb}[2]{\ifx!#1{\ref{nmb:#2}}\else\ifx|#1{#2}\else\ifx:#1{#2}\else{\label{nmb:#2}}\fi\fi\fi}
\allowdisplaybreaks

\ifcase\pdfoutput

\else
\usepackage[pdftex,raiselinks=true,backref=true
]{hyperref}%
\hypersetup{
bookmarksnumbered=true,%
bookmarksopen=true,%
bookmarksopenlevel=1
}
\fi

\providecommand{\mapsfrom}{\kern.2em%
\setbox0=\hbox{$\leftarrow$\kern-.10em\rule[0.26mm]{0.1mm}{1.3mm}}\box0%
\kern.3em}

\allowdisplaybreaks
\date{\today}

\title[Denjoy-Carleman differentiable perturbation]{Denjoy-Carleman differentiable perturbation of 
polynomials and unbounded operators
}

\author
[A.~Kriegl, P.W.~Michor, A.~Rainer]
{Andreas Kriegl, Peter W. Michor, and Armin Rainer}

\address{Andreas Kriegl: Fakult\"at f\"ur Mathematik, Universit\"at Wien, 
Nordbergstrasse~15, A-1090 Wien, Austria}
\email{andreas.kriegl@univie.ac.at}

\address{Peter W. Michor: Fakult\"at f\"ur Mathematik, Universit\"at Wien, 
Nordbergstrasse~15, A-1090 Wien, Austria}
\email{peter.michor@univie.ac.at}

\address{Armin Rainer: Department of Mathematics, University of Toronto, 
40 St.\ George Street, Toronto, Ontario, Canada M5S 2E4}

\email{armin.rainer@univie.ac.at}

\thanks{PM was supported by FWF-Project P~21030-N13.
AR was supported by FWF-Project J2771}


\begin{document}
\maketitle

\begin{abstract} 
Let $t\mapsto A(t)$ for $t\in T$ be a $C^M$-mapping with values unbounded 
operators with compact resolvents and common domain 
of definition which are self-adjoint or normal. 
Here $C^M$ stands for $C^\omega$ (real analytic), a quasianalytic or non-quasianalytic Denjoy-Carleman class, 
$C^\infty$, or a H\"older continuity class $C^{0,\alpha}$.
The parameter domain $T$ is either $\mathbb R$ or $\mathbb R^n$ or an infinite dimensional convenient vector space.
We prove and review results on $C^M$-dependence on $t$ of the eigenvalues and eigenvectors of $A(t)$.
\end{abstract}

\begin{proclaim}{Theorem}
Let $t\mapsto A(t)$ for $t\in T$ be a parameterized family of unbounded 
operators in a Hilbert space $H$ with common domain of definition and with 
compact resolvent.

If $t\in T=\mathbb R$ and all $A(t)$ are self-adjoint then the following holds:
\begin{enumerate}
\item[(A)] If $A(t)$ is real analytic in $t\in \mathbb R$, then 
the eigenvalues and the eigenvectors of $A(t)$ may be 
parameterized real analytically in $t$. 
\item[(B)] If $A(t)$ is quasianalytic of class $C^Q$ in $t\in \mathbb R$, then 
the eigenvalues and the eigenvectors of $A(t)$ may be 
parameterized $C^Q$ in $t$. 
\item[(C)] If $A(t)$ is non-quasianalytic of class $C^L$ in $t\in \mathbb R$ and if
no two unequal continuously parameterized 
eigenvalues meet of infinite order at any $t\in \mathbb R$, 
then the eigenvalues and the eigenvectors of $A(t)$ can be 
parameterized $C^L$ in $t$.
\item[(D)] If $A(t)$ is $C^\infty$ in $t\in \mathbb R$ and if
no two unequal continuously parameterized 
eigenvalues meet of infinite order at any $t\in \mathbb R$, 
then the eigenvalues and the eigenvectors of $A(t)$ can be 
parameterized $C^\infty$ in $t$.
\item[(E)] If $A(t)$ is $C^\infty$ in $t\in \mathbb R$, then 
the eigenvalues of $A(t)$ may be parameterized
twice differentiably in $t$. 
\item[(F)] If $A(t)$ is $C^{1,\alpha}$ in $t\in \mathbb R$ for some $\alpha>0$, 
then the eigenvalues of $A(t)$ may be 
parameterized in a $C^1$ way in $t$. 
\end{enumerate}
If $t\in T=\mathbb R$ and all $A(t)$ are normal then the following holds:
\begin{enumerate}
\item[(G)] If $A(t)$ is real analytic in $t\in \mathbb R$, then 
for each $t_0\in \mathbb R$ and for each eigenvalue $\lambda$ of $A(t_0)$ there exists 
$N\in \mathbb N$ such that the eigenvalues near $\lambda$ of $A(t_0\pm s^N)$ and their eigenvectors can be 
parameterized real analytically in $s$ near $s=0$.
\item[(H)] If $A(t)$ is $C^Q$ in $t\in \mathbb R$, then 
for each $t_0\in \mathbb R$ and for each eigenvalue $\lambda$ of $A(t_0)$ there exists 
$N\in \mathbb N$ such that the eigenvalues near $\lambda$ of 
$A(t_0\pm s^N)$ and their eigenvectors can be 
parameterized $C^Q$ in $s$ near $s=0$.
\item[(I)] If $A(t)$ is $C^L$ in $t\in \mathbb R$, then 
for each $t_0\in \mathbb R$ and for each eigenvalue $\lambda$ of $A(t_0)$ at which no two of the 
unequal continuously arranged eigenvalues (see \cite[II.5.2]{Kato76}) meet of infinite order,
there exists 
$N\in \mathbb N$ such that the eigenvalues near $\lambda$ of 
$A(t_0\pm s^N)$ and their eigenvectors can be 
parameterized $C^L$ in $s$ near $s=0$.
\item[(J)] If $A(t)$ is $C^\infty$ in $t\in \mathbb R$, then 
for each $t_0\in \mathbb R$ and for each eigenvalue $\lambda$ of $A(t_0)$ at which no two of the 
unequal continuously arranged eigenvalues (see \cite[II.5.2]{Kato76}) meet of infinite order,
there exists 
$N\in \mathbb N$ such that the eigenvalues near $\lambda$ of 
$A(t_0\pm s^N)$ and their eigenvectors can be 
parameterized $C^\infty$ in $s$ near $s=0$.
\item[(K)] If $A(t)$ is $C^\infty$ in $t\in \mathbb R$ and no two of the unequal continuously 
parameterized eigenvalues meet of infinite order at any $t\in \mathbb R$, then
the eigenvalues and the eigenvectors of $A(t)$ can be parameterized by absolutely continuous 
functions, locally in $t$. 
\end{enumerate}
If $t\in T=\mathbb R^n$ and all $A(t)$ are normal then the following holds:
\begin{enumerate}
\item[(L)] If $A(t)$ is $C^\omega$ or $C^Q$ in $t\in \mathbb R^n$, then
for each $t_0\in \mathbb R^n$ and for each eigenvalue $\lambda$ 
of $A(t_0)$, there exist
a finite covering $\{\pi_k : U_k \to W\}$ of a neighborhood $W$ of $t_0$, where each 
$\pi_k$ is a composite of finitely many mappings each of which is either a 
local blow-up along a $C^\omega$ or $C^Q$ submanifold or a local power substitution,
such that the eigenvalues and the eigenvectors of $A(\pi_k(s))$ can be chosen $C^\omega$ or $C^Q$ in $s$.
If $A$ is self-adjoint, then we do not need power substitutions. 
\item[(M)] If $A(t)$ is $C^\omega$ or $C^Q$ in $t\in \mathbb R^n$,
then the eigenvalues and their eigenvectors of $A(t)$ can be 
parameterized by functions which are special functions of bounded variation (SBV), see 
\cite{AmbrosioDeGiorgi88} or \cite{AFP00}, locally 
in $t$. 
\end{enumerate}
If $t\in T\subseteq E$, a $c^\infty$-open subset in an infinite dimensional convenient vector space then the
following holds:
\begin{enumerate}
\item[(N)] For $0<\alpha\le 1$, if $A(t)$ is $C^{0,\alpha}$ (H\"older continuous of exponent $\alpha$) in 
$t\in T$ and all $A(t)$ are self-adjoint, then the eigenvalues of $A(t)$ may be 
parameterized in a $C^{0,\alpha}$ way in $t$. 
\item[(O)] For $0<\alpha\le 1$, if $A(t)$ is $C^{0,\alpha}$ (H\"older continuous of exponent $\alpha$) in 
$t\in T$ and all $A(t)$ are normal, then we have:
For each $t_0\in T$ and each eigenvalue $z_0$ of $A(t_0)$ consider a simple closed 
$C^1$-curve $\gamma$ in the resolvent set of $A(t_0)$ enclosing only $z_0$ among all 
eigenvalues of $A(t_0)$. Then for $t$ near $t_0$ in the $c^\infty$-topology on $T$, 
no eigenvalue of $A(t)$ lies on $\gamma$.
Let $\lambda(t)=(\lambda_1(t),\dots,\lambda_N(t))$ be the $N$-tuple of all eigenvalues (repeated according 
to their multiplicity) of $A(t)$ inside of $\gamma$. 
Then $t\mapsto \lambda(t)$ is $C^{0,\alpha}$ for $t$ near $t_0$ with respect to the non-separating metric
\[
d(\lambda,\mu) = \min_{\sigma\in\mathcal S_N} \max_{1\le i\le N} |\lambda_i - \mu_{\sigma(i)}|
\]
on the space of $N$-tuples.
\end{enumerate}
\end{proclaim}

Part {\rm (A)} is due to Rellich \cite{Rellich42V} in 1942, 
see also \cite{Baumgaertel72} and \cite[VII,~3.9]{Kato76}.
Part {\rm (D)} has been proved in \cite[7.8]{AKLM98}, see 
also \cite[50.16]{KM97}, in 1997, which contains also a different proof 
of {\rm (A)}. {\rm (E)} and {\rm (F)} have been proved in \cite{KM03} in 
2003. 
(G) was proved in \cite[7.1]{RainerAC}; it can be proved as {\rm (H)} with some obvious 
changes, but it is not a special case since $C^\omega$ does not correspond to a sequence which is an 
$\mathcal L$-intersection (see \cite{KMRq}).
{\rm (J)} and {\rm (K)} were proved in \cite[7.1]{RainerAC}.
{\rm (N)} was proved in \cite{KMR}.

The purpose of this paper is to prove 
the remaining parts 
\thetag{B}, 
\thetag{C}, 
\thetag{H}, 
\thetag{I}, 
\thetag{L}, 
\thetag{M}, and 
\thetag{O}. 

\subsection*{Definitions and remarks }
Let $M=(M_k)_{k \in \mathbb{N}=\mathbb{N}_{\ge0}}$ be an increasing sequence ($M_{k+1}\ge M_k$) 
of positive real numbers with $M_0=1$.
Let $U \subseteq \mathbb{R}^n$ be open.
We denote by $C^M(U)$ the set of all $f \in C^\infty(U)$ such that, for each 
compact $K \subseteq U$,
there exist positive constants $C$ and $\rho$ such that
\begin{equation*}
|\partial^\alpha f(x)| \le C \, \rho^{|\alpha|} \, |\alpha|! \, M_{|\alpha|}
\quad\text{ for all }\alpha \in \mathbb{N}^n\text{ and }x \in K.
\end{equation*}
The set $C^M(U)$ is a \emph{Denjoy--Carleman class} of functions on $U$.
If $M_k=1$, for all $k$, then $C^M(U)$ coincides with the ring $C^\omega(U)$
of real analytic functions
on $U$. In general, $C^\omega(U) \subseteq C^M(U) \subseteq C^\infty(U)$.

Here $Q=(Q_k)_{k\in\mathbb N}$ is a sequence as above which is 
quasianalytic, log-convex, and which is also an $\mathcal L$-intersection,
see \cite{KMRq} or 
\cite{KMRc} and references therein. 
Moreover, $L=(L_k)_{k\in\mathbb N}$ is a sequence as above which is non-quasianalytic and log-convex. 

That $A(t)$ is a real analytic, $C^M$ (where $M$ is either $Q$ or $L$), $C^\infty$, or $C^{k,\alpha}$ family of 
unbounded operators means the following:
There is a dense subspace $V$ of the Hilbert space $H$ 
such that $V$ is the domain of definition of each $A(t)$, and such 
that $A(t)^*=A(t)$ in the self-adjoint case, or $A(t)$ has closed graph and 
$A(t)A(t)^*=A(t)^*A(t)$ wherever defined in the normal case. 
Moreover, we require 
that $t\mapsto \langle A(t)u,v\rangle$ 
is of the respective differentiability class for each $u\in V$ and $v\in H$. 
>From now on we treat only $C^M=C^\omega$, $C^M$ for $M=Q$, $M=L$, and $C^M=C^{0,\alpha}$. 

This implies that $t\mapsto A(t)u$ is of the same class $C^M(E,H)$ (where $E$ is either $\mathbb R$ 
or $\mathbb R^n$) or is in $C^{0,\alpha}(E,H)$ (if $E$ is a convenient vector space)
for each $u\in V$ 
by \cite[2.14.4, 10.3]{KM97} for $C^\omega$,
by \cite[3.1, 3.3, 3.5]{KMRc} for $M=L$,
by \cite[1.10, 2.1, 2.3]{KMRq} for $M=Q$,
and by \cite[2.3]{KM97}, \cite[2.6.2]{FK88} or \cite[4.14.4]{Faure91} 
for $C^{0,\alpha}$ because $C^{0,\alpha}$ can be described by boundedness conditions only 
and for these the uniform boundedness principle is valid.

A sequence of functions 
$\lambda_i$ is said to {\it parameterize the eigenvalues, if
for each $z\in \mathbb C$ 
the cardinality $|\{i: \lambda_i(t)=z\}|$ equals
the multiplicity of $z$ as eigenvalue of $A(t)$.} 

Let $X$ be a $C^\omega$ or $C^Q$ manifold.
A \emph{local blow-up $\Phi$} over an open subset $U$ of $X$ means the composition 
$\Phi = \iota \circ \varphi$ of a blow-up $\varphi : U' \to U$ with center a $C^\omega$ or $C^Q$ submanifold and 
of the inclusion $\iota : U \to X$.
A \emph{local power substitution} is a mapping $\Psi: V \to X$ of the form 
$\Psi = \iota \circ \psi$, where $\iota : W \to X$ is the inclusion of a coordinate chart $W$ of $X$ and 
$\psi : V \to W$ is given by 
\[
(y_1,\ldots,y_q) = ((-1)^{\epsilon_1} x_1^{\gamma_1},\ldots,(-1)^{\epsilon_q} x_q^{\gamma_q}),
\]
for some $\gamma=(\gamma_1,\ldots,\gamma_q) \in (\mathbb{N}_{>0})^q$ and all $\epsilon = (\epsilon_1,\ldots,\epsilon_q) \in \{0,1\}^q$, 
where $y_1,\ldots,y_q$ denote the coordinates of $W$ (and $q = \dim X$). 

This paper became possible only after some of the results of 
\cite{KMRc} and \cite{KMRq} were proved, 
in particular the uniform boundedness principles. The wish to prove the results of this 
paper was the main motivation for us to work on \cite{KMRc} and \cite{KMRq}.

\subsection*{Applications }
Let $X$ be a compact $C^Q$ manifold and let $t\mapsto g_t$ be a 
$C^Q$-curve of $C^Q$ Riemannian metrics on $X$. Then we get 
the corresponding $C^Q$ curve $t\mapsto \Delta(g_t)$
of Laplace-Beltrami operators on $L^2(X)$. By theorem (B) the
eigenvalues and eigenvectors can be arranged $C^Q$.
Question: Are the eigenfunctions then also $C^Q$? 

Let $\Omega$ be a bounded region in $\mathbb R^n$ with $C^Q$ boundary, and 
let $H(t)=-\Delta + V(t)$ be a $C^Q$-curve of Schr\"odinger 
operators with varying $C^Q$ potential and Dirichlet boundary conditions.
Then the eigenvalues and eigenvectors can be arranged $C^Q$.
Question: Are the eigenvectors viewed as eigenfunctions then also in $C^Q(\Omega\times \mathbb R)$? 

\subsection*{Example
} 
This is an elaboration of \cite[7.4]{AKLM98} and \cite[Example]{KM03}.
Let $S(2)$ be the vector space of all symmetric real 
$(2\times 2)$-matrices. 
We use the $C^L$-curve lemma \cite[3.6]{KMRc} or \cite[2.5]{KMRq}:
{\it
There exists a converging sequence of reals $t_n$ with the following 
property:
Let $A_n, B_n\in S(2)$ be any sequences which 
converge 
fast to 0, i.e., for each $k\in \mathbb N$ the sequences $n^kA_n$ and $n^kB_n$ are 
bounded in $S(2)$. 
Then there exists a curve $A\in C^L(\mathbb R,S(2))$ 
such that $A(t_n+s)=A_n+sB_n$ for 
$|s|\le \frac1{n^2}$, for all $n$.}

We use it for 
\begin{align*}
A_n := \frac 1{2^{n^2}} 
\begin{pmatrix}
1 & 0 \\
0 & -1\\
\end{pmatrix},\quad 
B_n := \frac 1{2^{n^2}\,s_n} 
\begin{pmatrix}
0 & 1 \\
1 & 0\\
\end{pmatrix},\quad 
\text{ where } s_n := 2^{n-n^2}\le \frac 1{n^2}.
\end{align*}
The eigenvalues of 
$A_n+tB_n$ and their derivatives are
\begin{equation*}
\lambda_n(t) = \pm\frac 1{2^{n^2}} \sqrt{1+(\tfrac t{s_n})^2},\quad
\lambda_n'(t) = \pm\frac {2^{n^2-2n}t}{\sqrt{1+(\frac t{s_n})^2}}.
\end{equation*}
Then
\begin{align*}
\frac{\lambda'(t_n+s_n)-\lambda'(t_n)}{s_n^\alpha} &=
\frac{\lambda_n'(s_n)-\lambda_n'(0)}{s_n^\alpha}
=\pm\frac {2^{n^2-2n}s_n}{s_n^\alpha\sqrt{2}}\\
&=\pm\frac{2^{n(\alpha(n-1)-1)}}{\sqrt{2}} \to \infty \text{ for }\alpha>0.
\end{align*}
So condition (in \thetag{C}, \thetag{D}, \thetag{I}, \thetag{J}, and \thetag{K}) 
that no two unequal continuously parameterized eigenvalues meet of infinite order 
cannot be dropped.
By \cite[2.1]{AKLM98}, we may always find a 
twice differentiable square root of a non-negative smooth function, 
so that the eigenvalues $\lambda$ are functions 
which are twice differentiable but not $C^{1,\alpha}$ for any $\alpha>0$. 

Note that the normed eigenvectors cannot be chosen continuously in 
this example (see also example \cite[\S 2]{Rellich37I}). 
Namely, we have
\begin{equation*}
A(t_n)=A_n=\frac1{2^{n^2}}\begin{pmatrix} 1 & 0 \\ 0 &-1 \end{pmatrix},\qquad
A(t_n+s_n)=A_n +s_n\,B_n=\frac1{2^{n^2}}\begin{pmatrix} 1 & 1 \\ 1 &-1 \end{pmatrix}.
\end{equation*}

\begin{proclaim}{Resolvent Lemma}
Let $C^M$ be any of $C^\omega$, $C^Q$, $C^L$, $C^\infty$, or $C^{0,\alpha}$, and let $A(t)$ be normal. 
If $A$ is $C^M$ then the resolvent 
$(t,z)\mapsto (A(t)-z)^{-1}\in L(H,H)$ is $C^M$ on its natural domain,
the global resolvent set
\[
\{(t,z)\in T\times\mathbb C: (A(t)-z):V\to H \text{ is invertible}\}
\]
which is open (and even connected).
\end{proclaim}

\begin{demo}{\bf Proof}
By definition the function $t\mapsto \langle A(t)v,u \rangle$ is of 
class $C^M$ for each $v\in V$ and $u\in H$. 
We may conclude that
the mapping $t\mapsto A(t)v$ is of class $C^M$ into $H$ 
as follows:
For $C^M=C^\infty$ we use \cite[2.14.4]{KM97}.
For $C^M=C^\omega$ we use in addition \cite[10.3]{KM97}. 
For $C^M=C^Q$ or $C^M=C^L$ we use \cite[2.1]{KMRq} and/or \cite[3.3]{KMRc} 
where we replace $\mathbb R$ by $\mathbb R^n$.
For $C^M=C^{0,\alpha}$ we use 
\cite[2.3]{KM97}, \cite[2.6.2]{FK88}, or \cite[4.1.14]{Faure91} 
because $C^{0,\alpha}$ can be described by
boundedness conditions only and for these the uniform boundedness
principle is valid.

For each $t$ consider the norm $\|u\|_t^2:=\|u\|^2+\|A(t)u\|^2$ on 
$V$.
Since $A(t)$ is closed, $(V,\|\quad\|_t)$ is again a 
Hilbert space with inner product 
$\langle u,v\rangle_t:=\langle u,v\rangle+\langle A(t)u,A(t)v\rangle$. 

{\it {\rm (1)} Claim (see \cite[in the proof of 7.8]{AKLM98}, 
\cite[in the proof of 50.16]{KM97}, or \cite[Claim 1]{KM03}). 
All these norms $\|\quad\|_t$ on $V$ are equivalent, 
locally uniformly in $t$.
We then equip $V$ with one of the 
equivalent Hilbert norms, say $\|\quad\|_0$.}

We reduce this to $C^{0,\alpha}$.
Namely, note first that $A(t):(V,\|\quad\|_s)\to H$ is bounded since the 
graph of $A(t)$ is closed in $H\times H$, contained in $V\times H$ and thus 
also closed in $(V,\|\quad\|_s)\times H$. 
For fixed $u,v\in V$, the function
$t\mapsto \langle u,v\rangle_t=\langle u,v \rangle+\langle A(t)u,A(t)v 
\rangle$ is $C^{0,\alpha}$ since $t\mapsto A(t)u$ is it.
By the multilinear uniform boundedness principle 
(\cite[5.18]{KM97} or \cite[3.7.4]{FK88}) the 
mapping $t\mapsto \langle \quad,\quad\rangle_t$ is $C^{0,\alpha}$ into the 
space of bounded sesquilinear forms on $(V,\|\quad\|_s)$ for each fixed 
$s$. 
Thus the inverse image of $\langle \quad,\quad \rangle_s + \frac12(\text{unit ball})$ in 
$L(\overline{(V,\|\quad\|_s)} \oplus (V,\|\quad\|_s);\mathbb C)$
is a $c^\infty$-open neighborhood $U$ of $s$ in $T$.
Thus $\sqrt{1/2}\|u\|_s\le \|u\|_t\le \sqrt{3/2}\|u\|_s$ for all $t\in U$, i.e.\  
all Hilbert norms $\|\quad\|_t$ are locally uniformly 
equivalent, and claim {\rm (1)} follows. 

By the 
linear uniform boundedness theorem 
we see that $t\mapsto A(t)$ is in $C^M(T, L(V,H))$ as follows (here it suffices to use 
a set of linear functionals which together recognize bounded sets instead of the whole dual):
For $C^M=C^\infty$ we use \cite[1.7 and 2.14.3]{KM97}.
For $C^M=C^\omega$ we use in addition \cite[9.4]{KM97}. 
For $C^M=C^Q$ or $C^M=C^L$ we use \cite[2.2 and 2.3]{KMRq} and/or \cite[3.5]{KMRc} 
where we replace $\mathbb R$ by $\mathbb R^n$.
For $C^M=C^{0,\alpha}$ see above.

If for some $(t,z)\in T\times\mathbb C$ the bounded operator 
$A(t)-z:V\to H$ is invertible, then this is true locally with respect to the $c^\infty$-topology on 
the product which is the product topology by \cite[4.16]{KM97}, and 
$(t,z)\mapsto (A(t)-z)^{-1}:H\to V$ is $C^M$, by the chain rule,
since inversion is real analytic on the Banach space $L(V,H)$.
\qed\end{demo}

Note that $(A(t)-z)^{-1}:H\to H$ is a compact operator 
for some (equivalently any) $(t,z)$ if and only if the inclusion 
$i:V\to H$ is compact, since $i=(A(t)-z)^{-1}\circ(A(t)-z): V\to H\to H$.

\begin{proclaim}{Polynomial proposition}
Let $P$ be a curve of polynomials 
\[
P(t)(x)=x^n-a_1(t)x^{n-1}+\dots+(-1)^na_n(t),\quad t\in \mathbb R.
\]
\begin{enumerate}
\item[(a)] If $P$ is hyperbolic (all roots real) and if 
the coefficient functions $a_i$ are all $C^Q$
then there exist $C^Q$ functions $\lambda_i$ which parameterize all roots. 
\item[(b)] If $P$ is hyperbolic (all roots real), 
if the coefficient functions $a_i$ are $C^L$ and no two of the different roots meet of 
infinite order, then there exist $C^L$ functions 
$\lambda_i$ which parameterize all roots. 
\item[(c)] If the coefficient functions $a_i$ are $C^Q$, 
then for each $t_0$ there exists $N\in \mathbb N$ such that the roots of 
$s\mapsto P(t_0\pm s^{N})$ can be parameterized $C^Q$ in $s$ for $s$ near 0.
\item[(d)] If the coefficient functions $a_i$ are $C^L$ and no two of the different roots meet of 
infinite order, 
then for each $t_0$ there exists $N\in \mathbb N$ such that the roots of 
$s\mapsto P(t_0\pm s^{N})$ can be parameterized $C^L$ in $s$ for $s$ near 0.
\end{enumerate}
All $C^Q$ or $C^L$ solutions differ by permutations.
\end{proclaim}

The proof of parts \thetag{a} and \thetag{b} is exactly as in \cite{AKLM98} where the corresponding results were 
proven for $C^\infty$ instead of $C^L$, and for $C^\omega$ instead of $C^Q$. For this we need only 
the following properties of $C^Q$ and $C^L$:
\begin{itemize}
\item They allow for the implicit function theorem (for \cite[3.3]{AKLM98}).
\item They contain $C^\omega$ and are closed under composition (for \cite[3.4]{AKLM98}).
\item They are derivation closed (for \cite[3.7]{AKLM98}).
\end{itemize}
Part \thetag{a} is also in \cite[7.6]{CC04} which follows \cite{AKLM98}. It also follows from 
the multidimensional version \cite[6.10]{RainerQA} since blow-ups in dimension 1 are trivial.
The proofs of parts \thetag{c} and \thetag{d} are exactly as in \cite[3.2]{RainerAC} where the 
corresponding result was proven for $C^\omega$ instead of $C^Q$, and for
$C^\infty$ instead of $C^L$, if none of the different roots meet of infinite 
order. For these we need the properties of $C^Q$ and $C^L$ listed above. 

\begin{proclaim}{Matrix proposition}
Let $A(t)$ for $t\in T$ be a family of $(N\times N)$-matrices.
\begin{enumerate}
\item[(e)] If $T=\mathbb R\ni t\mapsto A(t)$ is a $C^Q$-curve of Hermitian matrices,
then the eigenvalues and the eigenvectors can be chosen $C^Q$.
\item[(f)] If $T=\mathbb R\ni t\mapsto A(t)$ is a $C^L$-curve of Hermitian matrices
such that no two eigenvalues meet of infinite order,
then the eigenvalues and the eigenvectors can be chosen $C^L$.
\item[(g)] If $T=\mathbb R\ni t\mapsto A(t)$ is a $C^L$-curve of normal matrices
such that no two eigenvalues meet of infinite order,
then for each $t_0$ there exists $N_1\in \mathbb N$ such that the eigenvalues and eigenvectors of 
$s\mapsto A(t_0\pm s^{N_1})$ can be parameterized $C^L$ in $s$ for $s$ near 0.
\item[(h)] Let $T\subseteq\mathbb R^n$ be open and let $T\ni t\mapsto A(t)$ be a $C^\omega$ or $C^Q$-mapping of 
normal matrices.
Let $K \subseteq T$ be compact. 
Then there exist a neighborhood $W$ of $K$, and
a finite covering $\{\pi_k : U_k \to W\}$ of $W$, where each 
$\pi_k$ is a composite of finitely many mappings each of which is either a 
local blow-up along a $C^\omega$ or $C^Q$ submanifold or a local power substitution,
such that the eigenvalues and the eigenvectors of $A(\pi_k(s))$ can be chosen $C^\omega$ or $C^Q$ in $s$.
Consequently, the eigenvalues and eigenvectors of $A(t)$ are locally special functions of bounded 
variation (SBV).
If $A$ is a family of Hermitian matrices, then we do not need power substitutions. 
\end{enumerate}
\end{proclaim}

The proof of the matrix proposition in case \thetag{e} and \thetag{f} is exactly as in 
\cite[7.6]{AKLM98}, using the polynomial proposition and properties of $C^Q$ and $C^L$.
Item \thetag{g} is exactly as in \cite[6.2]{RainerAC}, using the polynomial proposition and properties of $C^L$.
Item \thetag{h} is proved in \cite[9.1 and 9.6]{RainerQA}, see also \cite{KurdykaPaunescu08}.

\begin{demo}{\bf Proof of the theorem}
We have to prove parts 
\thetag{B}, 
\thetag{C}, 
\thetag{H}, 
\thetag{I}, 
\thetag{L}, 
\thetag{M}, and 
\thetag{O}. 
So let $C^M$ be any of $C^\omega$, $C^Q$, $C^L$, or $C^{0,\alpha}$, and let $A(t)$ be normal.
Let $z$ be an eigenvalue of $A(t_0)$ of multiplicity $N$. 
We choose a simple closed $C^1$ curve $\gamma$ in the resolvent set of 
$A(t_0)$ for fixed $t_0$ enclosing only $z$ among all eigenvalues of $A(t_0)$. 
Since the global resolvent set is open, see the resolvent lemma, no eigenvalue 
of $A(t)$ lies on $\gamma$, for $t$ near $t_0$.
By the resolvent lemma, $A: T\to L((V,\|\quad\|_0),H)$ is $C^M$, thus also
\begin{equation*}
t\mapsto -\frac1{2\pi i}\int_\gamma (A(t)-z)^{-1}\;dz =: P(t,\gamma) = P(t)
\end{equation*}
is a $C^M$ mapping. Each $P(t)$ is a projection, namely onto the direct sum of all 
eigenspaces corresponding to eigenvalues of $A(t)$ in the interior of $\gamma$, with finite rank.
Thus the rank must be constant: It is easy to see 
that the 
(finite) rank cannot fall locally, and it cannot increase, since the 
distance in $L(H,H)$ of $P(t)$ to the subset of operators of 
rank $\le N=\operatorname{rank}(P(t_0))$ is continuous in $t$ and is either 
0 or 1. 

So for $t$ in a neighborhood $U$ of $t_0$
there are equally many eigenvalues in the 
interior of $\gamma$, and we may call them 
$\lambda_i(t)$ for $1\le i\le N$ (repeated 
with multiplicity).

Now we consider the family of $N$-dimensional complex vector spaces 
$t\mapsto P(t)(H)\subseteq H$, for $t\in U$. They form a $C^M$ Hermitian vector subbundle over $U$ 
of $U\times H\to U$:
For given $t$, choose $v_1,\dots v_N\in H$ such that the $P(t)v_i$ are linearly independent and thus span $P(t)H$. 
This remains true locally in $t$. Now we use the Gram Schmidt orthonormalization 
procedure (which is $C^\omega$) for the $P(t)v_i$ to obtain a local orthonormal $C^M$ frame of the bundle. 

Now $A(t)$ maps $P(t)H$ to itself; in a $C^M$ local frame it is given by a normal $(N\times N)$-matrix 
parameterized $C^M$ by $t\in U$. 

Now all local assertions of the theorem follow:
\begin{itemize}
\item[\thetag{B}]
Use the matrix proposition, part \thetag{e}.
\item[\thetag{C}]
Use the matrix proposition, part \thetag{f}.
\item[\thetag{H}]
Use the matrix proposition, part \thetag{h}, and note that in dimension 1 blowups are trivial.
\item[\thetag{I}]
Use the matrix proposition, part \thetag{g}.
\item[\thetag{L,M}]
Use the matrix proposition, part \thetag{h}, for $\mathbb R^n$.
\item[\thetag{O}]
We use the following \newline
{\bf Result} {\it \emph{(\cite{BhatiaDavisMcIntosh83}, \cite[VII.4.1]{Bhatia97})}
Let $A,B$ be normal $(N \times N)$-matrices and
let $\lambda_i(A)$ and 
$\lambda_i(B)$ for $i=1,\dots, N$ denote the respective 
eigenvalues.
Then 
\[
\min_{\sigma\in\mathcal S_N}\max_j |\lambda_j(A)-\lambda_{\sigma(j)}(B)| \le C \|A-B\|
\]
for a universal constant $C$ with $1 < C < 3$.
Here $\|\quad\|$ is the operator norm.}
\end{itemize}

Finally, it remains to extend the local choices to global ones for the cases \thetag{B} and 
\thetag{C} only. There $t\mapsto A(t)$ is $C^Q$ or $C^L$, respectively, which imply both 
$C^\infty$, and no two different eigenvalues meet of infinite order. So we may apply 
\cite[7.8]{AKLM98} (in fact we need only the end of the proof) to conclude that 
the eigenvalues can be chosen $C^\infty$ on $T=\mathbb R$, uniquely up to a global 
permutation. By the local result above they are then $C^Q$ or $C^L$. The same proof then gives us, 
for each eigenvalue $\lambda_i:T \to \mathbb R$ with generic multiplicity $N$, a unique $N$-dimensional 
smooth vector subbundle of $\mathbb R\times H$ whose fiber over $t$ consists of eigenvectors for the 
eigenvalue $\lambda_i(t)$. In fact this vector bundle is $C^Q$ or $C^L$ by the local result above, 
namely the matrix proposition, part \thetag{e} or \thetag{f}, respectively.
\qed\end{demo}

\def\cprime{$'$}
\providecommand{\bysame}{\leavevmode\hbox to3em{\hrulefill}\thinspace}
\providecommand{\MR}{\relax\ifhmode\unskip\space\fi MR }
\providecommand{\MRhref}[2]{%
  \href{http://www.ams.org/mathscinet-getitem?mr=#1}{#2}
}
\providecommand{\href}[2]{#2}

\end{document}